\documentclass[10pt]{amsart}

\issueinfo{00}% volume number
  {}%         % issue number
  {}%         % month
  {2005}%     % year

\copyrightinfo{1995}  
{American Mathematical Society}% copyright holder

\usepackage{amsfonts,amssymb}
\usepackage{amscd}
\usepackage[all]{xy}

\newtheorem{theorem}{Theorem}[section]
\newtheorem{lemma}[theorem]{Lemma}
\newtheorem{Prop}[theorem]{Proposition}
\newtheorem{corollary}[theorem]{Corollary}

\theoremstyle{definition}
\newtheorem{definition}[theorem]{Definition}

\theoremstyle{remark}
\newtheorem{remark}[theorem]{Remark}

\numberwithin{equation}{section}

\def\dsp{\displaystyle}
\newcommand{\ncmd}{\newcommand}
\ncmd{\beq}{\begin{equation}}
\ncmd{\eeq}{\end{equation}}
\ncmd{\beqa}{\begin{eqnarray}}
\ncmd{\eeqa}{\end{eqnarray}}
\ncmd{\beqan}{\begin{eqnarray*}}
\ncmd{\eeqan}{\end{eqnarray*}}
\def\demo{{\bf Proof: }}
\def\qed{$\hfill\Box$}

\begin{document}

\setlength{\parindent}{0pt}

\title[Automorphisms]{Automorphisms of $\mathbb{C}^k$
and associated compact complex manifolds.}

%    Information for first author
\author{RENAUD Julie}
%    Address of record for the research reported here
\address{Centre de Mathematiques et d'Informatique\\
 39 rue Joliot-Curie. 13453 Marseille Cedex 13}
%    Current address
%\curraddr{Department of Mathematics and Statistics,
%Case Western Reserve University, Cleveland, Ohio 43403}
\email{renaud@cmi.univ-mrs.fr, jujuren@hotmail.com}


\date{July 7, 2005.}

%\dedicatory{This paper is dedicated to our advisors.}


\begin{abstract} In this paper, we first construct $k$-dimensional compact
complex manifolds from automorphisms of
$\mathbb{C}^k$ which admit a fixed attracting point at infinity. Then,
we charactize the fundamental group as well as the universal covering
of the attracting basin of this fixed point thanks to a generalization of
the method described by T. Bousch in his thesis \cite{B}.
\end{abstract}

\medskip

\maketitle

%\tableofcontents

\section{Introduction}

\medskip

In the paper \cite{DO}, the authors give the construction of surfaces
which contain a global spherical shell (G.S.S) from the attracting basin
of an H\'enon automorphism that is an automorphism such that:
$$\begin{array}{cccc}
H\ :\ &\mathbb{C}^2 &\longrightarrow & \ \mathbb{C}^2 \\
& \left(\begin{array}{l} x\\ y\\ \end{array}\right)
&\longmapsto & \left(\begin{array}{l} x^2+c-ay\\ \qquad x\\
\end{array}\right)
\end{array}$$

The properties that the surfaces satisfy are consequences of results
obtained by J. Hubbard and R.Oberste-Vorth in
\cite{HOV}. In this paper, the authors are interested in the
dynamical system constituted by an H\'enon mapping and the set
$\mathbb{C}^2 \backslash K^+$ where $$K^+ := \Big\{(x,y)\mid H^{\circ
n}(x,y)\not\rightarrow \infty \ {\rm when} \ n\rightarrow \pm
\infty \Big\}.$$
They give a topologic model as the projective limit of tori, then
infer a holomorphic model. In his thesis, T. Bousch
\cite{B} considers the same dynamical system but tackles the problem in a
different way. He finds the same topologic model as a consequence of
various calculations among others, asymptotic developments. This method
valid in dimension $2$ also holds in dimension $3$
as shown by the example of the automorphism $H(x,\ y,\ z) = (y,\ z,\ z^2
+ x)$ given in his thesis.\\

Our aim is on the one hand to generalize to dimension $k$ the previous
construction, from $H$ a polynomial automorphism of $\mathbb{C}^k$
of degree $d$ which, once extended to $\mathbb{P}^{k}$ given in
homogeneous coordinates $[z_1:\ z_2:\
\cdots:\ z_k:\ t]$,
contracts the hyperplane at infinity $\{t = 0\}$ minus the
indeterminacy set onto a fixed attracting point $p$, in order to get a
manifold with a global spherical shell, what we do in the first part; on
the other hand to establish in a second part the structure of the
attracting basin by using asymptotic developments.
\smallskip

In these two parts, we will study the class of automorphisms given
by:

$$ H(z_1,\ z_2,\ \cdots,\ z_k) = \left\{ \begin{array}{l}
z_1^d + \alpha_2 z_2 + \cdots\ + \alpha_k z_k\\
z_3 \\
\ \vdots\\
\ \vdots\\
z_k\\
z_1
\end{array}\right.$$
with $\alpha_2\neq 0,\
\vert\alpha_2\vert < d $ and $\alpha_3,\ \cdots,\ \alpha_k$
any constant in $\mathbb{C}$.\\
Let $U^+$ be the attracting basin at infinity of the automorphism we
consider, that is the set:

$$U^+ := \Big\{w\in\mathbb{C}^k, \ \displaystyle\lim_{n\to\infty} H^n(w) =
p
\Big\}.$$
\smallskip
The result we get is summed up in the following theorem:\\

{\bf Theorem:}
$$\pi_1\ ( U^{+}) \simeq \mathbb{Z}\ \Big[\ \frac{1}{d}\ \Big]$$

$${\rm and \ } \ \widetilde{U^{ +}} = \mathbb{H}\ \times \mathbb{C}^{
k-1}$$

{\it where $\mathbb{H}$ denotes the Poincar\'e half-plane.}\\

\smallskip

In the third part, we will deal with the particular case of
quadratic automorphisms of $\mathbb{C}^{3}$ which have a fixed
attracting point at infinity. The frame of this work is the
classification of all these automorphisms into five classes established
by Fornaess and Wu \cite{FW}. We give the list of all the quadratic
automorphisms of $\mathbb{C}^{3}$ from which we can construct a threefold
with a global spherical shell.

\smallskip

\section{Construction of a manifold with a global spherical
shell.}

\medskip

We can do in dimension $k$ the analogue of the construction of surfaces
with global spherical shell as done in \cite{D} and proceed as follows.\\

We extend the automorphism $H$ to $G$ birational, from
$\mathbb{P}^k(\mathbb{C})$ to
$\mathbb{P}^k(\mathbb{C})$ given once again in homogeneous coordinates
by $[z_1:\ \cdots:\ z_k:\ t]$.

$$\begin{array}{cccc}
G : &\mathbb{P}^k(\mathbb{C}) &\longrightarrow &
\mathbb{P}^k(\mathbb{C}) \\
&\ &\ &\ \\
&[z_1: \cdots: z_k: t] &\longmapsto &\big[z_1^d +
\displaystyle\sum_{j=2}^{k} \alpha_j z_j t^{d-1} : z_3
t^{d-1}: \cdots:
z_k t^{d-1}: z_1 t^{d-1}: t^d\big].\\
\end{array}$$

Denote $I = \big\{z_1 = t = 0\big\}$ the set of indeterminacy;
the hyperplane $\big\{ t = 0\big\}$ minus $I$ is
flattened onto the fixed point $p = \big[1:0:\cdots:0\big]$
which is superattracting.

Indeed, the birational mapping $G$ written in coordinates in $U_1 =
\{z_1\neq 0\}$, becomes

$$F(\zeta_2,\cdots,\zeta_{k+1}) =
\left(\displaystyle\frac{\zeta_3\zeta_{k+1}^{d-1}}{D},\cdots,
\displaystyle\frac{\zeta_k\zeta_{k+1}^{d-1}}{D},
\displaystyle\frac{\zeta_{k+1}^{d-1}}{D},
\displaystyle\frac{\zeta_{k+1}^{d}}{D}\right)$$\\

with $D := 1 + \alpha_2\
\zeta_2\zeta_{k+1}^{d-1} +
\cdots +
\alpha_k\ \zeta_k\zeta_{k+1}^{d-1}$ and where\\ $
\big(\zeta_2,\cdots,\zeta_{k+1}\big) = \Big(\dsp\frac{z_2}{z_1},
\cdots,\dsp\frac{z_k}{z_1},\dsp\frac{t}{z_1} \Big) = (0,\cdots,0)$
corresponds to the point $p$.\\

The Jacobian matrix of $F$ at $(0,\ \cdots,\ 0)$ is

$$J =
\left(\begin{array}{c c c}
0 & \cdots & 0 \\ &
\\ \vdots & 0 & \vdots \\ &
\\ 0 &\cdots & 0\\
\end{array}\right)$$\\

therefore the point $p$ is superattracting.\\

The map $F$ also satisfies the following:

\begin{lemma} The germ $F: (\mathbb{C}^k,0) \longrightarrow
(\mathbb{C}^k,0)$ can be written as $\pi\eta$ with $\pi$ a composition of
several blow-ups above the unit ball $B$ and
$\eta$ a local biholomorphism.
\end{lemma}

\demo Let $\eta$ be the application
$$\eta : (\zeta_{2},\ \cdots,\ \zeta_k)
\longmapsto
\Big (\zeta_3,\ \cdots,\ \zeta_k,\displaystyle\frac{-\alpha_2\zeta_2 -
\cdots - \alpha_k\zeta_k}{D},\
\zeta_{k+1}\Big) $$ and let $\pi_j$ be the $2d -1$ blow-ups
defined in local coordinates by:
$$\pi_1 : (u_1,\ \cdots,\ u_k) \longmapsto (u_1 u_{k-1},\ \cdots,\
u_{k-2} u_{k-1},\ u_{k-1},\ u_k u_{k-1})$$

$$ \pi_2 = \cdots = \pi_d: (u_1,\ \cdots,\ u_k) \longmapsto (u_1,\
\cdots,\ u_{k-2},\ u_{k-1} u_{k},\ u_{k})$$

$$\qquad \pi_{d+1} : (u_1,\ \cdots,\ u_k) \longmapsto
(u_1,\
\cdots,\ u_{k-2},\ u_{k-1} u_{k} + 1,\ u_{k}) \qquad \qquad $$

$$ \pi_{d+2} = \cdots = \pi_{2d-1}: (u_1,\ \cdots,\ u_k) \longmapsto
(u_1,\
\cdots,\ u_{k-2},\ u_{k-1} u_{k},\ u_{k}).$$

Let $\pi$ be the composition $\pi = \pi_1\circ\pi_2\circ
\cdots\circ\pi_{2d-1}$. Then the germ $F$ can be written as $F = \pi\eta$
and it is easy to check that $\eta$ is a local biholomorphism.\qed\\

For a small $\varepsilon > 0$, let $B$ be $B_{\varepsilon} :=
\left\{ w \in \mathbb{C}^k \mid \vert \vert w \vert\vert < 1 + \varepsilon
\right\}$.\\

The map $\eta\pi: \pi^{-1}(B_{\varepsilon}) \longrightarrow
\pi^{-1}(B_{\varepsilon})$ sends biholomorphically a neighbourhood of
the boundary $\partial\pi^{-1}(B)$ of $\pi^{-1}(B)$ in $B_{\varepsilon}$,
onto a neighbourhood of
the boundary $\partial(\eta(B))$ of $\eta(B)$; by gluing up with $\eta\pi$
these holomorphic neighbourhoods, we get a compact complex manifold of
dimension $k$ called $X$, with a global spherical shell.

\smallskip

\begin{remark} Let $\Gamma$ be the maximal divisor of $X$, then
$X\backslash\Gamma$ is isomorphic to the quotient $U^+/\langle H\rangle$.

\end{remark}

>From now on, we try to understand the structure of $U^+$,
that is the fundamental group $\pi_1(U^+)$ and the universal covering
$\widetilde{U^+}$ of the attracting basin of the automorphism $H$.

\medskip

\section{Structure of the attracting basin.}

\subsection{Fundamental group.}$\ $\\

\smallskip

First, note that $U^{+}$ can be written
in the form:
$$U^{+} = \displaystyle\bigcup_n \ H^{-n}(V^{+}) \ {\rm where} \
V^{+} = \Big\{(z_1,\cdots,\ z_k)\in \mathbb{C}^{k} \mid \vert z_1\vert
 	 	> Max\{R,\ \vert
z_2\vert\ ,\ \cdots,\ \vert
z_k\vert\}\Big\},$$ where $R >> 0$.

We begin by constructing a holomorphic function which will be usefull in
what follows:

\begin{Prop} There exists a holomorphic function $\varphi : V^+
\longrightarrow
\mathbb{C}\backslash \overline{\Delta}$, with $\Delta$ the unit disk,
such that:
\begin{enumerate}
\item For all $(z_1,\cdots,z_k) \in V^+$,
$\varphi\Big(H\big(z_1,\cdots,z_k \big)\Big) =
\Big(\varphi(z_1,\cdots,z_k)\Big)^d. \qquad (*)$\\

\item $\varphi(z_1,\cdots,z_k) \sim z_1$ when $\vert z_1\vert
\longrightarrow +\infty$ in $V^+. \hfill (**)$
\end{enumerate}

\end{Prop}

\demo For a point $z=(z_1,\cdots, z_k)$, we write $H^{\circ
n}(z)$ in the form $$H^{\circ n}(z) = (f_{1,n}(z),\ \cdots,\
f_{k,n}(z)), \ {\rm for} \ n \in\mathbb{Z}.$$

We will show that, for each $n$, we can find a $d^n$th root
called $\varphi_n(z)$ of $f_{1,n}(z)$ on $V^+$ and the sequence
$\varphi_n(z)$ converges uniformly to the required holomorphic function
$\varphi$. From the definition of $H$, we have:

$$f_{1,1} = z_1^d + \alpha_2 z_2 + \cdots + \alpha_k z_k \ {\rm \ or}$$

$$\displaystyle\frac{f_{1,1}}{z_1^d} = 1 + \alpha_2\displaystyle\frac{
z_2}{z_1^d} + \cdots +
\alpha_k\displaystyle\frac{ z_k}{z_1^d} \ .$$

If $z \in V^{+}$ , then $\vert f_{1,1}\vert > R$, hence $1 + \alpha_2
\displaystyle\frac{z_2}{z_1^d} + \cdots +
\alpha_k\displaystyle\frac{ z_k}{z_1^d}$ admits a logarithm $\beta(z)$ and
$f_{1,1} = z_1^d\exp\Big(\beta(z)\Big)$.\\

The logarithm is of the form: $\beta(z) = \ (\alpha_2 z_2 + \cdots +
\alpha_k z_k) O\Big(\displaystyle\frac{1}{z_1^d}\Big)$ and $\beta$ est
bounded on $V^{+}$ because

\begin{eqnarray*} \displaystyle\vert
\displaystyle\frac{\alpha_2 z_2 + \cdots + \alpha_k z_k}{z_1^d}\vert &<&
\displaystyle\frac{\vert \alpha_2\ \vert \vert z_2\vert + \cdots + \vert
\alpha_k\
\vert \vert z_k\vert }{\vert z_1\vert^d} \\
&\ &\\
&<& \displaystyle\frac{\vert
\alpha_2\ \vert + \cdots + \vert
\alpha_k\
\vert }{\vert R\vert^{d-1}}.\\
\end{eqnarray*}

In the above equation, we replace $z$ by $H^{\circ n-1}(z) =
\big(f_{1,n-1}(z),\ \cdots, f_{k,n-1}(z)\big)$, and we obtain

$$f_{1,n}(z) =
f_{1,1}\Big(H^{\circ (n-1)}(z)\Big)= f_{1,n-1}^{d}(z)\ \exp\Big(
\beta\big( H^{\circ (n-1)} (z)\big)\Big).$$

Using similar relations for $H^{\circ n-2}, H^{\circ n-3}, \cdots, $ we
have

$$f_{1,n}(z) = z_1^{d^{ n}} exp \left\{ \ d^{ n-1}\
\beta (z) \ + \ \dots\ + \ \beta \big(H^{\circ (n-1)}(z)\big)\right\}.$$

The $d^{ n}$-th root appears to be the function:

$$\varphi_n := z_1 \ exp \ \left\{ \frac{1}{d} \ \beta(z) \ + \dots\ +
\frac{1}{d^{ n}}\ \beta\big(H^{\circ (n-1)}(z)\big)\right\}$$

function for which the series converges normally since
$\beta$ is bounded on $V^{+}$. The function $\varphi(z) :=
\displaystyle\lim_{n\to \infty} \varphi_{n}(z)$ satisfies the two
expected properties since $$\varphi_n\big(H(z)\big) =
\big(\varphi_{n+1}(z)\big)^d.$$
We can determine $\varphi$ more precisely:

\begin{eqnarray*} \varphi(w) &=& z_1 \ exp\left(\alpha_2
\displaystyle\frac{z_2}{z_1^d} + \cdots +
\alpha_k\displaystyle\frac{ z_k}{z_1^d} \right) + \cdots\\
&=& z_1\ \left(\ 1 +
\alpha_2\ O\Big(\displaystyle\frac{z_2}{z_1^d}\Big) + \cdots + \alpha_k\
O\Big(\displaystyle\frac{z_k}{z_1^d}\Big)\right)\\
&=& z_1 \ + \alpha_2\
O\Big(\displaystyle\frac{z_2}{z_1^{d-1}}\Big) + \cdots + \alpha_k\
O\Big(\displaystyle\frac{z_k}{z_1^{d-1}}\Big).\qquad
\qquad \qquad \hfill \Box
\end{eqnarray*}

\begin{corollary}\label{co2} Let $(z_{0} ,\ \cdots, \ z_{-k+1})$ be an
element of $ \ V^{+}$,  $  U := \big(\varphi(z_{0},\ \cdots, \
z_{-k+1})\big)^{d^n}$ and let the sequences $(z_{1,n})_n, \cdots,
(z_{k,n})_n$ be defined by
$$H^{\circ n}\big(z_{0} ,\ \cdots, z_{-k+1}\big) =: \big(z_{1,n} ,\
\cdots, z_{k,n}\big).$$
The sequence $(z_{1,n})_n$ admits as a first estimation the following
development:

$$ z_{1,n} = U + O\ \Big(U^{d^{-k+j-1} + 1 - d}\Big) \quad {\rm if} \
\alpha_{j+1} = \cdots = \alpha_k = 0 \ {\rm and} \ \alpha_j\neq 0.$$
\end{corollary}

\demo Let $(z_{0} ,\ \cdots,\ z_{-k+1})$ be in $V^{+}$ and let us
consider:
$$ \ H^{\circ n}\big(z_{0} ,\ \cdots, \ z_{-k+1}\big) = \big(z_{1,n}
,\
\cdots, z_{k,n}\big)$$ we have
$$\varphi(z_{1,n} ,\
\cdots, z_{k,n}) = \varphi\big(H^{\circ n}(z_{0} ,\ \cdots,\
z_{-k+1})\big)\stackrel{(*)}{=} \big(\varphi(z_{0} ,\ \cdots,\
z_{-k+1})\big)^{d^{
n}} = U .$$

An equivalent of $z_{1,n}$ is $U$ thanks to property $(*)$
of the previous proposition. The first estimation is obtained thanks
to the function $\varphi(z)$.\qed

\smallskip

\begin{Prop} There exists a closed $1$-form $\omega$ on
$U^{+}$ such that $$H^{*}\omega = d\ \omega.$$
\end{Prop}

\demo We define $\omega := {\rm d}(\log \varphi)$ or
$\displaystyle\frac{{\rm d}\varphi}{\varphi}$ on $V^{+}$. We have :

\begin{eqnarray*} \varphi\ (H(z)) &=& \varphi^{d}(z)\\
{\rm or\  } \ \ H^{*}\varphi
&=& \varphi^{d}\ \ \qquad \quad \ \
(1)\\ {\rm by\ differentiating\ } \ \ H^{*}{\rm d}\varphi
&=& d \ \varphi^{d-1}\ {\rm d}\varphi\
\qquad (2)\\ {\rm by\ dividing\ } (2) \ {\rm by\ } (1) \ \ \
H^{*}\displaystyle\frac{{\rm d}\varphi}{\varphi}
&=& d\ \displaystyle\frac{{\rm
d}\varphi}{\varphi}\\
{\rm thus}\ \ H^{*}\omega
&=& d \ \omega \ {\rm on} \
V^{+}.
\end{eqnarray*}
In general, we define $\omega:= \dsp\frac{1}{d^n}
(H^n)^*\omega$ sur $H^{-n}(V^+)$.\qed\\

This closed form $\omega$ allows us to study from now on the fundamental
group $\pi_1(U^+)$ of $U^+$.

\medskip
For a closed curve $C$ in $U^{+}$, we set

$$ \alpha (C) := \displaystyle\frac{1}{2 i \pi}\displaystyle\int_C
\omega.$$

Since $\omega$ is a closed $1-$form, the number $\alpha(C)$ is
determined by the homotopy class of $C$.

\smallskip
\begin{Prop} For any closed curve $C$ in $U^{+}$, the following
assertions holds:

\begin{enumerate}

\item $\ \ \alpha(H(C)) = d\ \alpha(C)$,

\item $\ \ \alpha(C)$ belongs to $\mathbb{Z}
\Big[\displaystyle\frac{1}{d}
\Big]$,

\item \ $C$ is zero-homotopic in $U^{+}$ if and only if
$\alpha(C) = 0$,

\item for any $ r \in \mathbb{Z}\ \Big[\displaystyle\frac{1}{d}\Big],\
$ there exists a closed curve $C$ in $U^{+}$ such that $$\alpha ( C ) =
r.$$
\end{enumerate}
\end{Prop}

\demo\\

%\begin{enumerate}
$1.$ $\alpha(H(C)) = \displaystyle\frac{1}{2 i
\pi}\displaystyle\int_{H(C)} \omega =
\displaystyle\frac{1}{2 i \pi}\displaystyle\int_C
H^{*}\omega \stackrel{(*)}{=} \displaystyle\frac{1}{2 i \pi}\ d
\displaystyle\int_C
\omega = d \ \alpha(C)$.\\

$\ $

$2.$ Assume $C$ is in $V^{+}$.\\
In $V^{+},$ $\varphi$
is in the form $\varphi = z_1\ e^{ \gamma(z)}$ hence
$\displaystyle\frac{{\rm d} \varphi}{\varphi} = \omega =
\displaystyle\frac{{\rm d} z_1}{z_1} \ + \ {\rm d} \gamma$.\\
For a closed curve $C$ in $V^{+}$,

$$\alpha(C) =
\displaystyle\frac{1}{2 i
\pi}\displaystyle\int_{C} \omega = \displaystyle\frac{1}{2 i
\pi}\displaystyle\int_{C} \displaystyle\frac{{\rm d} z_1}{z_1} +
\displaystyle\frac{1}{2 i
\pi}\displaystyle\int_{C} {\rm d} \gamma$$

since $C$ is closed, $\displaystyle\int_C {\rm d} \gamma = 0$, and \
$\displaystyle\frac{1}{2 i
\pi}\displaystyle\int_{C} \omega = \displaystyle\frac{1}{2 i
\pi}\displaystyle\int_{C} \displaystyle\frac{{\rm d} z_1}{z_1}$.\\

This vanishes if $C$ is zero-homotopic in $V^{+}$.\\

Let $C$ denote a closed  curve in $U^{+}$. There exists an integer $n$
such that
$H^{n} (C) \subset V^{+}$ that is such that $C \subset H^{-n}(V^{+})$. We
know that $\alpha(H^{n} (C))$ is an integer hence

$$\alpha (C ) = \displaystyle\frac{1}{d^{ n}} \ \alpha( H^{n}
(C))\ \ {\rm belongs \ to\  } \mathbb{Z}\
\Big[\displaystyle\frac{1}{d}\Big].$$
If $\alpha (C)= 0 $, then $\alpha (H^{n}(C)) = 0$ and hence $H^{n}(C)$
is zero-homotopic in $V^{+}$.\\ Since $H^n : H^{-n} (V^{+})
\longrightarrow V^{+} $ is a biholomorphism, $C$ is zero-homotopic in
$H^{-n} (V^{+})$.
Let $C_0$ denote the closed curve in $V^{+}$ given
by:

$$\begin{array}{cccc}
C_0\ :\ &[\ 0\ ,\ 1\ ] &\longrightarrow & \ V^{+}\\
&\ &\ &\ \\
&t &\longmapsto & \big( 2 R e^{2 i \pi
t}\ ,\ 0\ ,\cdots,\ 0\big)
\end{array}$$
we have
$$\alpha\ ( C_0 ) = \displaystyle\frac{1}{2 i
\pi}\displaystyle\int_{C_0} \omega = 1 \ \ {\rm et\ } \alpha\ ( m\ C_0 )
= m.$$
Let $r:= \displaystyle\frac{m}{d^{ n}}\ \in \mathbb{Z}\
\Big[\displaystyle\frac{1}{d}\Big]$, we choose a closed curve $C:=
H^{ - n}(m\ C_0)$; it is such that
$$\alpha(C) = \alpha\ \big(H^{ - n}(m\ C_0)\big) =
\displaystyle\frac{1}{d^{ n}} \ \alpha\ (m\ C_0) =
\displaystyle\frac{m}{d^{ n}} = r.$$
This proves the proposition.\qed\\
%\end{enumerate}

The proposition shows that $\alpha$ is an isomorphism of $\pi_1 (
U^{+} )$ onto $ \mathbb{Z}\ \big[\displaystyle\frac{1}{d}\Big]$ and the
fundamental group is determined.

\smallskip

\begin{remark} Note that we could have used the Green functions
to find $\varphi$, but the estimation of $z_{1,n}$ would not have been
precise enough for our needs in the following paragraph.
\end{remark}

\smallskip

\subsection{Fundamental group of $U^+$.}$\ $\\

\medskip

>From now on, we are looking for the universal covering of $U^{+}$.
We apply the method of the asymptotic development introduced by
T.Bousch in his thesis. We have to distinguish two cases according the
degree $d$ equals $2$ or is greater than $3$.\\
We show the following proposition:\\

\begin{Prop}\label{P} There exists one and only one holomorphic
function
$$G\ : \ (\mathbb{C} \backslash \overline{\Delta}) \times \mathbb{C}^{k-1}
\longrightarrow \mathbb{C}^{k} \backslash K^+ = U^+$$ for which the
following properties hold:\\

$1.$ $G$ is a locally trivial bundle with discrete countable fibers.\\

$2.$ Let $\omega :\ (\mathbb{C} \backslash \overline{\Delta})
\times\
\mathbb{C}^{k-1}
\longrightarrow (\mathbb{C} \backslash \overline{\Delta}) \times \
\mathbb{C}^{k-1}$ be defined by

\begin{eqnarray*}
\omega(v,s_1,\cdots,s_{k-1}) &=& \Bigg(v^{ d}, \Lambda_1 \Big(s_1+
\dsp\frac{v^{(d^{k}-1)\cdot d^{k-2}}}{k-1}-
\displaystyle\sum_{l=3}^{k}\displaystyle\frac{\alpha_l}{d(k-1)}
v^{\big(-d^{l-3}+d^{l+k-4}+(d-1)d^{k-2}\big)}\Big)\\
&\ &\\
& ,\cdots, & \Lambda_{k-1} \Big(s_{k-1} +
\dsp\frac{v^{(d^{k}-1)\cdot d^{k-2}}}{k-1}-
\displaystyle\sum_{l=3}^{k}\displaystyle\frac{\alpha_l}{d(k-1)}
v^{\big(-d^{l-3}+d^{l+k-4}+(d-1)d^{k-2}\big)} \Big)\Bigg)
\end{eqnarray*}
(where $\Lambda_1,\cdots,\Lambda_{k-1}$ are solutions of
$dr^{k-1} + \alpha_2 = 0$), in the case $d\geq 3$;\\

or, if $d=2$, be defined by

\begin{eqnarray*}
\omega(v,s_1,\cdots,s_{k-1}) &=& \Bigg(v^{ 2}, \Lambda_1 \Big(s_1+
\dsp\frac{v^{(2^{k}-1)\cdot 2^{k-1}}}{k-1}-
\displaystyle\sum_{l=3}^{k}\displaystyle\frac{\alpha_l}{2(k-1)}
v^{\big(-2^{l-2}+2^{l+k-3}+2^{k-1}\big)}\\
&\ &\\
&+& \displaystyle\frac{\alpha_3\alpha_k}{4(k-1)} v\Big)\ ,\cdots,
\Lambda_{k-1} \Big(s_{k-1} +
\dsp\frac{v^{(2^{k}-1)\cdot 2^{k-1}}}{k-1}\\
&\ &\\
&-&
\displaystyle\sum_{l=3}^{k}\displaystyle\frac{\alpha_l}{2(k-1)}
v^{\big(-2^{l-2}+2^{l+k-3}+2^{k-1}\big)} +
\displaystyle\frac{\alpha_3\alpha_k}{4(k-1)} v \Big)\Bigg)
\end{eqnarray*}
(where $\Lambda_1,\cdots,\Lambda_{k-1}$ are solutions of
$2 \cdot r^{k-1} + \alpha_2 = 0$)\\

then the following diagram is commutative:
\large{\[
\begin{CD}
\big(\mathbb{C}\ \backslash\ \overline{\Delta}\big)\ \times\
\mathbb{C}^{k-1} @>^{\rm{\omega}}
>>
 \big(\mathbb{C}\ \backslash\ \overline{\Delta}\big)\
\times\
\mathbb{C}^{k-1}\\
@VGVV @VVGV\\
\mathbb{C}^k\ \backslash\ K^+ = U^+ @>H>> \mathbb{C}^k\ \backslash\
K^+ = U^+
\end{CD}
\]}

\end{Prop}

Let us set, for $n \in \mathbb{Z}$, $(z_{1,n}, z_{2,n},\cdots, z_{k,n})
:= H^{\circ n} (z_0, z_{-1},\cdots, z_{-k+1})$.
To study the sequence of the iterates of $(z_0, z_{-1},\cdots, z_{-k+1})$
under
$H$, it is sufficient to study the sequence $(z_{1,n})_n$ which satisfies
an induction formula of the type:
$$ z_{1,n+1} = z_{1,n}^{d} + \alpha_2 z_{1,n-k+1} + \alpha_3
z_{1,n-k+2} + \cdots + \alpha_k z_{1,n-1}.$$

The idea of the demonstration of the proposition lies on the asymptotic
development of $z_{1,n}$. In what follows, the calculations are purely
formal. We will deal with the convergence later. We suppose that
$\alpha_k \neq 0$, if not, we replace $k$ by
$k-1$ for instance.\\

We begin with a preliminary lemma:\\

\begin{lemma}\label{L} $z_{1,n}$ admits the following asymptotic development
for
$d\geq 3$:

\beqan
z_{1,n} &=& U^{\displaystyle\frac{(d-1)d^{k-1}}{1 - d^{k-1}}} \ \Bigg[
\displaystyle\sum_{j=0}^{\infty} \ \Big(-\frac{\alpha_2}{d}\Big)^{ j} \
U^{
\displaystyle\frac{1-d^k}{1-d^{k-1}} \cdot (d^{1-k})^j}\\
&\ &\\
&\ & - \displaystyle\sum_{l=3}^{k}\displaystyle\frac{\alpha_l}{d}
\displaystyle\sum_{j=0}^{\infty} \ \Big(-\frac{\alpha_2}{d}\Big)^{ j} \
U^{
\displaystyle\frac{d^{l-k-1}-d^{l-2}-(d-1)}{1-d^{k-1}} \cdot
(d^{1-k})^j}\Bigg] + \cdots
\eeqan

and for $d=2$:
\beqan
z_{1,n} &=& U^{\displaystyle\frac{2^{k-1}}{1 - 2^{k-1}}} \ \Bigg[
\displaystyle\sum_{j=0}^{\infty} \ \Big(-\frac{\alpha_2}{2}\Big)^{ j} \
U^{
\displaystyle\frac{1-2^k}{1-2^{k-1}} \cdot (2^{1-k})^j}\\
&\ &\\
&- & \displaystyle\sum_{l=3}^{k}\displaystyle\frac{\alpha_l}{2}
\displaystyle\sum_{j=0}^{\infty} \ \Big(-\frac{\alpha_2}{2}\Big)^{ j} \
U^{
\displaystyle\frac{2^{l-k-1}-2^{l-2}-1}{1-2^{k-1}} \cdot
(2^{1-k})^j}\\
&\ &\\
&+ & \dsp\frac{\alpha_3\alpha_k}{4} \dsp\sum_{j=0}^{\infty} \
\Big(-\frac{\alpha_2}{2}\Big)^{ j} \ U^{
\displaystyle\frac{-1}{(1-2^{k-1})\cdot 2^{k-1}} \cdot
(2^{1-k})^j}\Bigg] + \cdots\\
\eeqan

\end{lemma}

\demo We prove the lemma for $d$ greater than $3$.\\
To establish the asymptotic development of $z_{1,n}$ thanks to the first
estimation, we isolate $z_{1,n}$ in the induction formula then we inject
the developments of $z_{1,n+1}$, $z_{1,n-1}$, $\cdots$, $z_{1,n-k+1}$ in
the relation which gives $z_{1,n}$; we suppose $\alpha_k \neq 0$.\\

We have: $$z_{1,n} = \big(\ z_{1,n+1} - \alpha_2 z_{1,n-k+1}
- \cdots - \alpha_k z_{1,n-1}\big)^{\frac{1}{d}}.$$

We refine the first estimation of the asymptotic development given by the
corollary. We begin with the estimation of $z_{1,n}$ and we use the fact
that from $n$ to $n+1$, $U$ is changed into $U^d$:

\begin{eqnarray*} z_{1,n} &=& U + O\big(U^{\frac{1}{d}-d+1} \big)\\
\ &\ & \\
{\rm donc\ }\ \ z_{1,n+1} &=& U^{ d} + O \big(U^{1-d^2+d} \big)\\
\ &\ & \\
\ z_{1,n-1} &=& U^{\frac{1}{d}} + O\big( U^{\frac{1-d^2+d}{d^2}}
\big)\\
\vdots &\ & \qquad\quad\vdots\\
\ &\ & \\
\vdots &\ & \qquad\quad\vdots\\
\ &\ & \\
z_{1,n-k+2} &=& U^{\frac{1}{d^{k-2}}} + O\big( U^{\frac{1-d^2+d}{d^{k-1}}}
\big)\\
\ &\ & \\
z_{1,n-k+1} &=& U^{\frac{1}{d^{k-1}}} + O\big( U^{\frac{1-d^2+d}{d^{k}}}
\big)\\
\end{eqnarray*}

thus

\begin{eqnarray*}
z_{1,n} &=& \left(U^{d} - \alpha_2 \ U^{\frac{1}{d^{k-1}}} -
\alpha_3
\ U^{\frac{1}{d^{k-2}}} - \cdots - \alpha_k \ U^{\frac{1}{d}} + O \Big( \
U^{\frac{1-d^2+d}{d^{k}}} \Big)\right)^{\frac{1}{d}}\\
&=& U\ \left(1 - \alpha_2 \
U^{\frac{1-d^k}{d^{k-1}}} - \alpha_3
\ U^{\frac{1-d^{k-1}}{d^{k-2}}} - \cdots\right.\\
&\ & \left.\qquad \qquad \qquad \qquad - \alpha_k \
U^{\frac{1-d^2}{d}} + O
\Big(
\ U^{\frac{1-d^2+d-d^{k+1}}{d^{k}}} \Big)\right)^{\frac{1}{d}}\\
&=& U\ \left(1 - \frac{\alpha_2}{d} \
U^{\frac{1-d^k}{d^{k-1}}} - \frac{\alpha_3}{d}
\ U^{\frac{1-d^{k-1}}{d^{k-2}}} - \cdots\right.\\
&\ & \left. \qquad \qquad \qquad \qquad \ - \frac{\alpha_k}{d} \
U^{\frac{1-d^2}{d}} + O
\Big(
\ U^{\frac{1-d^2+d-d^{k+1}}{d^{k}}} \Big) \right)\\
\end{eqnarray*}
\begin{eqnarray*}
\  &=& U - \frac{\alpha_2}{d} \
U^{\frac{1-d^k+d^{k-1}}{d^{k-1}}} - \frac{\alpha_3}{d}
\ U^{\frac{1-d^{k-1}+d^{k-2}}{d^{k-2}}} - \cdots\\
&\ & \qquad \qquad \qquad \qquad \ \ \ -
\frac{\alpha_k}{d}
\ U^{\frac{1-d^2+d}{d}} + O
\Big(
\ U^{\frac{1-d^2+d-d^{k+1}+d^k}{d^{k}}} \Big)\\
&=& U - \frac{\alpha_2}{d} \
U^{\frac{1-d^{k-1}(d-1)}{d^{k-1}}} - \frac{\alpha_3}{d}
\ U^{\frac{1-d^{k-2}(d-1)}{d^{k-2}}} - \cdots\\
&\ & \qquad \qquad \qquad \qquad \quad \ - \frac{\alpha_k}{d} \
U^{\frac{1-d(d-1)}{d}} + O
\Big(
\ U^{\frac{1-d(d-1)-d^k(d-1)}{d^{k}}} \Big)\\
\end{eqnarray*}

We use this new development to refine a little bit more and we
get $z_{1,n}$ in the form of a series (we are not confronted with the
problem of convergence since the calculations are formal):

\begin{eqnarray*}
z_{1,n} &=& U^{\displaystyle\frac{(d-1)d^{k-1}}{1 - d^{k-1}}} \ \Bigg[
\displaystyle\sum_{j=0}^{\infty} \ \Big(-\frac{\alpha_2}{d}\Big)^{ j} \
U^{
\displaystyle\frac{1-d^k}{1-d^{k-1}} \cdot (d^{1-k})^j}\\
&\ &\\
&\ & - \displaystyle\sum_{l=3}^{k}\displaystyle\frac{\alpha_l}{d}
\displaystyle\sum_{j=0}^{\infty} \ \Big(-\frac{\alpha_2}{d}\Big)^{ j} \
U^{
\displaystyle\frac{d^{l-k-1}-d^{l-2}-(d-1)}{1-d^{k-1}} \cdot
(d^{1-k})^j}\Bigg] + \cdots
\end{eqnarray*}
For $d=2$, a new term appears, the one associated with the coefficient
$\alpha_k$ in the development of $z_{1,n-k+2}$,
this is the reason why a term with a coefficient equal
to $\alpha_3\alpha_k$ appears.\qed

%\begin{re} This development is not entirely satisfactory.
%Tout
%d'abord, les exposants qui entrent en jeu dans le d\'eveloppement sont de
%la forme $p/2^k$ avec $k$ grand, ce qui pose le probl\`eme de la
%d\'efinition de $U^{ 2^{-k}}$. Ce probl\`eme est lev\'e si l'on
%consid\`ere que $U$ est un \'el\'ement non pas de
%$\mathbb{C}$ mais de $\widehat{\mathbb{C}}$ o\`u $\widehat{\mathbb{C}}$
%est la limite projective de $\mathbb{C}$ sous l'action de $z\ \mapsto
%z^2$. De plus, les exposants ne tendent pas vers $- \infty$ comme un vrai
%d\'eveloppement asymptotique mais s'accumulent vers $-2$.

\begin{remark}The development of $z_{1,n}$ depends not only on $U$. Moreover,
several sequences $(z_{1,n})$ which have the same $U$, differ by an error
which is about $$U^{ \displaystyle\frac{(d-1)d^{k-1}}{1
- d^{k-1}}}.$$
\end{remark}

To estimate the error, we write once again $z_{1,n} $:

\begin{eqnarray*}
z_{1,n} &=& U^{\displaystyle\frac{(d-1)d^{k-1}}{1 - d^{k-1}}} \ \Bigg[
\displaystyle\sum_{j=0}^{\infty} \ \Big(-\frac{\alpha_2}{d}\Big)^{ j} \
U^{
\displaystyle\frac{1-d^k}{1-d^{k-1}} \cdot (d^{1-k})^j}\\
&\ &\\
&\ & - \displaystyle\sum_{l=3}^{k}\displaystyle\frac{\alpha_l}{d}
\displaystyle\sum_{j=0}^{\infty} \ \Big(-\frac{\alpha_2}{d}\Big)^{ j} \
U^{
\displaystyle\frac{d^{l-k-1}-d^{l-2}-(d-1)}{1-d^{k-1}} \cdot
(d^{1-k})^j}\Bigg]\\ &\ &\\
&\ & \qquad \qquad \qquad \qquad + \ U^{\displaystyle\frac{(d-1)d^{k-1}}{1
- d^{k-1}}} \cdot v_n \ .
\end{eqnarray*}

We take once more the induction formula which gives $z_{1,n}$ and we
compare the coefficients.\\

Up to $\dsp\frac{d-1}{1-d^{k-1}}$, the coefficients exactly
compensate. When one compares the coefficients
of $U^{\frac{d-1}{1-d^{k-1}}}$, one obtains the induction formula:

$$ d \ v_n + \alpha_2 \ v_{n-k+1} = 0.$$

The sequences which satisfy this relation form a vector space
of dimension $k-1$, generated by the solutions of
$$d r^{k-1} +
\alpha_2 = 0,$$ a basis of which is noted
$\Lambda_1,\
\cdots,\
\Lambda_{k-1}$. We set:
\begin{eqnarray*}
z_{1,n} &=& U^{\displaystyle\frac{(d-1)d^{k-1}}{1 - d^{k-1}}} \ \Bigg[
\displaystyle\sum_{j=0}^{\infty} \ \Big(-\frac{\alpha_2}{d}\Big)^{ j} \
U^{
\displaystyle\frac{1-d^k}{1-d^{k-1}} \cdot (d^{1-k})^j}\\
&\ &\\
&- & \displaystyle\sum_{l=3}^{k}\displaystyle\frac{\alpha_l}{d}
\displaystyle\sum_{j=0}^{\infty} \ \Big(-\frac{\alpha_2}{d}\Big)^{ j} \
U^{
\displaystyle\frac{d^{l-k-1}-d^{l-2}-(d-1)}{1-d^{k-1}} \cdot
(d^{1-k})^j}\Bigg]\\
&\ &\\
&+ & \ U^{\displaystyle\frac{(d-1)d^{k-1}}{1
- d^{k-1}}} \cdot \Big(Q_1 + \cdots + Q_{k-1}\Big) \ \ {\rm where} \
\ Q_j := \big(\Lambda_j\big)^n\ q_{j,0}.
\end{eqnarray*}
For $d=2$, the induction formula which gives the error is $$2 v_n
+
\alpha_2 v_{n-k+1} = 0. $$
The sequences which satisfy this relation form a $(k-1)$-dimensionnal
vector space generated by the solutions of $2 r^{k-1} +
\alpha_2 = 0$, we denote $\Lambda_1,\ \cdots,\
\Lambda_{k-1}$ a basis. We set

\medskip

\beqan
z_{1,n} &=& U^{\displaystyle\frac{2^{k-1}}{1 - 2^{k-1}}} \ \Bigg[
\displaystyle\sum_{j=0}^{\infty} \ \Big(-\frac{\alpha_2}{2}\Big)^{ j} \
U^{
\displaystyle\frac{1-2^k}{1-2^{k-1}} \cdot (2^{1-k})^j}\\
&\ &\\
&- & \displaystyle\sum_{l=3}^{k}\displaystyle\frac{\alpha_l}{2}
\displaystyle\sum_{j=0}^{\infty} \ \Big(-\frac{\alpha_2}{2}\Big)^{ j} \
U^{
\displaystyle\frac{2^{l-k-1}-2^{l-2}-1}{1-2^{k-1}} \cdot
(2^{1-k})^j}\\
&\ &\\
&+ & \dsp\frac{\alpha_3\alpha_k}{4} \dsp\sum_{j=0}^{\infty} \
\Big(-\frac{\alpha_2}{2}\Big)^{ j} \ U^{
\displaystyle\frac{-1}{(1-2^{k-1})\cdot 2^{k-1}} \cdot
(2^{1-k})^j}\Bigg]\\
&\ &\\
& +& \ U^{\displaystyle\frac{2^{k-1}}{1
- 2^{k-1}}} \cdot \Big(Q_1 + \cdots + Q_{k-1}\Big) \ \ {\rm where} \
\ Q_j := \big(\Lambda_j\big)^n\ q_{j,0}.\\
\eeqan

Once we have established the developement of $z_{1,n}$ as well as the
sequences $(Q_1,\cdots,Q_{k-1})$, we need to determinate the universal
covering of $U^{+}$. That's what we do from now on
through the resolution of a functional equation. We can already state
that $\ \widetilde{U^{+}} =
\mathbb{H}\times\mathbb{C}^{k-1}$. \\

We have just seen that $z_{1,n}$ only depends on $k$
parameters $U,\ Q_1,\ \cdots,\ Q_{k-1}$ which become, when one passes
from $n$ to $n+1$, \ $(U^d,\ \Lambda_1 Q_1,\ \cdots,\ \Lambda_{k-1}
Q_{k-1})$. This transformation will be used as a model and  helps
us to find a map:

$$\begin{array}{cccc}
f :\ &\big(\mathbb{C} \backslash \overline{\Delta}\big) \times
\mathbb{C}^{k-1} &\longrightarrow & \mathbb{C} \\
&\ &\ &\ \\
& (u,\ q_1,\ \cdots,\ q_{k-1}) &\longmapsto & f(u,\
\ q_1,\ \cdots,\ q_{k-1})\\
\end{array}$$

such that the sequence $\Big(f\big(u^{d^{n}} ,\ \Lambda_1^n
q_1,\cdots,\Lambda^n_{k-1} q_{k-1}\big)\Big)_n$ satisfies the relation:
$$\forall \ (u,\ q_1,\ \cdots,\ q_{k-1}) \in
\big(\mathbb{C}
\backslash
\overline{\Delta}\big) \times \mathbb{C}^{k-1},$$
\begin{eqnarray*}
f\big(u^d,\Lambda_1 q_1, \cdots,\Lambda_{k-1} q_{k-1}\big) &=&
\Big(f\big(u, q_1, \cdots, q_{k-1}\big)\Big)^d + \alpha_2
f\Bigg(u^{\frac{1}{d^{k-1}}}, \frac{q_1}{\Lambda_1^{k-1}}, \cdots,
\frac{q_{k-1}}{\Lambda_{k-1}^{k-1}} \Bigg)\\
&\ & \qquad + \ \cdots\ + \alpha_k
f\Bigg(u^{\frac{1}{d}}, \frac{q_1}{\Lambda_1}, \cdots,
\frac{q_{k-1}}{\Lambda_{k-1}} \Bigg).
\end{eqnarray*}

We begin with

\begin{eqnarray*}
f_0\big(u, q_1, \cdots, q_{k-1}\big)&=&
u^{\displaystyle\frac{(d-1)d^{k-1}}{1 - d^{k-1}}} \ \Bigg[
\displaystyle\sum_{j=0}^{\infty} \ \Big(-\frac{\alpha_2}{d}\Big)^{ j} \
u^{
\displaystyle\frac{1-d^k}{1-d^{k-1}} \cdot (d^{1-k})^j}\\
&\ &\\
&- & \displaystyle\sum_{l=3}^{k}\displaystyle\frac{\alpha_l}{d}
\displaystyle\sum_{j=0}^{\infty} \Big(-\frac{\alpha_2}{d}\Big)^{ j} \
u^{
\dsp\frac{d^{l-k-1}-d^{l-2}-(d-1)}{1-d^{k-1}} \cdot (d^{1-k})^j}\Bigg]\\
&\ &\\
&+ & \ u^{\displaystyle\frac{(d-1)d^{k-1}}{1
- d^{k-1}}} \cdot \Big(q_1 + \cdots + q_{k-1}\Big)
\end{eqnarray*}

and we make several approximations by rewriting the functional equation
as follows $\Big($ this new writing makes the fractional exponents, which
are not well defined on $\big(\mathbb{C} \backslash
\overline{\Delta}\big),$ disappear $\Big)$:

\begin{eqnarray*}
f\big(u^{d^k},\Lambda_1^{k} q_1, \cdots,\Lambda^{k}_{k-1}
q_{k-1}\big) &=&
\Big(f\big(u^{d^{k-1}}, \Lambda_1^{k-1}q_1, \cdots,
\Lambda_{k-1}^{k-1}q_{k-1}\big)\Big)^d\\
&+ & \alpha_2
f\Big(u, q_1,\cdots, q_{k-1} \Big)
+ \cdots \\
&+ & \alpha_k
f\Big(u^{d^{k-2}}, \Lambda_{1}^{k-2} q_1, \cdots,
\Lambda_{k-1}^{k-2} q_{k-1} \Big).
\end{eqnarray*}

We set

$$\begin{array}{cccc}
\omega : &\big(\mathbb{C}\backslash \overline{\Delta}\big)
\times\mathbb{C}^{k-1} &\longrightarrow &\big(\mathbb{C}\backslash
\overline{\Delta}\big)
\times\mathbb{C}^{k-1}
\\
&\ &\ &\\
&(u, q_1,\cdots, q_{k-1}) &\longmapsto &\big(u^d,\Lambda_1 q_1,
\cdots,\Lambda_{k-1} q_{k-1}\big),
\end{array}$$

and we are led to study a simpler functional equation:

$$f = \dsp\frac{1}{\alpha_2} \Big[ f\circ\omega^{\circ k} -
\big( f\circ\omega^{\circ (k-1)}\big)^d - \alpha_3 f\circ \omega -
\cdots - \alpha_k f\circ \omega^{\circ (k-2)}\Big].$$

We get a formal series $f(u,q_1,\cdots,q_{k-1})$ which satisfies
this functional equation, but we cannot cope with the problem of
convergence. The formula giving $f_0$ with
$u\in
\big(\mathbb{C}\backslash \overline{\Delta}\big)$ et $(q_1,\cdots,q_{k-1})
\in
\mathbb{C}^{k-1}$ has sense as a formal series but converges only
if $\vert\alpha_2\vert <d$.\\

Moreover, in order to obtain a series whose exponents are integer, we
gather the fractional exponents by setting, in the case $d\geq 3$:

\begin{eqnarray*}
s_1 &=& q_1 +
\dsp\frac{1}{k-1}
\displaystyle\sum_{j=1}^{\infty} \Lambda_1^{ j} \
u^{
\displaystyle\frac{1-d^k}{1-d^{k-1}} \cdot d^{-j}}\\
&\ &\\
&\ & \quad -
\dsp\frac{1}{k-1}\displaystyle\sum_{l=3}^{k}\displaystyle\frac{\alpha_l}{d}
\displaystyle\sum_{j=1}^{\infty} \Lambda_1^{ j} \
u^{
\dsp\frac{d^{l-k-1}-d^{l-2}-(d-1)}{1-d^{k-1}} \cdot d^{-j}}\\
\vdots &\ & \qquad \qquad \qquad \vdots \qquad \qquad \qquad \vdots\\
\vdots &\ & \qquad \qquad \qquad \vdots \qquad \qquad \qquad \vdots\\
s_{k-1} &=& q_{k-1} +
\dsp\frac{1}{k-1}
\displaystyle\sum_{j=1}^{\infty} \Lambda_{k-1}^{ j} \
u^{
\displaystyle\frac{1-d^k}{1-d^{k-1}} \cdot d^{-j}}\\
&\ &\\
&\ &\qquad -
\dsp\frac{1}{k-1}\displaystyle\sum_{l=3}^{k}\displaystyle\frac{\alpha_l}{d}
\displaystyle\sum_{j=1}^{\infty} \Lambda_{k-1}^{ j} \
u^{
\dsp\frac{d^{l-k-1}-d^{l-2}-(d-1)}{1-d^{k-1}} \cdot d^{-j}}
\end{eqnarray*}

We put the previous expressions of $q_1,\cdots,q_{k-1}$ into $z_{1,n}$,
this allows us to gather the fractional exponents.

\begin{eqnarray*}
z_{1,n} &=&
u^{\displaystyle\frac{(d-1)d^{k-1}}{1 - d^{k-1}}} \
\Bigg[ u^{\displaystyle\frac{1-d^{k}}{1 - d^{k-1}}} -
\displaystyle\sum_{l=3}^{k}\displaystyle\frac{\alpha_l}{d}
\ u^{
\dsp\frac{d^{l-k-1}-d^{l-2}-(d-1)}{1-d^{k-1}}}\\
&\ & +\
\displaystyle\sum_{j=1}^{\infty} \ \Big(-\frac{\alpha_2}{d}\Big)^{ j} \
u^{
\displaystyle\frac{1-d^k}{1-d^{k-1}} \cdot (d^{1-k})^j}\\
&\ &\\
&\ & - \ \displaystyle\sum_{l=3}^{k}\displaystyle\frac{\alpha_l}{d}
\displaystyle\sum_{j=1}^{\infty} \Big(-\frac{\alpha_2}{d}\Big)^{ j} \
u^{
\dsp\frac{d^{l-k-1}-d^{l-2}-(d-1)}{1-d^{k-1}} \cdot (d^{1-k})^j}\\
&\ &\\
&\ & + \ s_1 + \cdots + s_{k-1} -
\dsp\frac{1}{k-1}
\displaystyle\sum_{j=1}^{\infty} \Big(\Lambda_{1}^{ j} +
\cdots + \Lambda_{k-1}^{ j}\Big) \ u^{
\displaystyle\frac{1-d^k}{1-d^{k-1}} \cdot d^{-j}}\\
&\ & +\
\dsp\frac{1}{k-1}\displaystyle\sum_{l=3}^{k}\displaystyle\frac{\alpha_l}{d}
\displaystyle\sum_{j=1}^{\infty} \Big(\Lambda_{1}^{ j} +
\cdots + \Lambda_{k-1}^{ j}\Big) \
u^{
\dsp\frac{d^{l-k-1}-d^{l-2}-(d-1)}{1-d^{k-1}} \cdot
d^{-j}}\Bigg]
\end{eqnarray*}

\begin{eqnarray*}
z_{1,n} &=& u^{\displaystyle\frac{(d-1)d^{k-1}}{1 - d^{k-1}}} \
\Bigg[ u^{\displaystyle\frac{1-d^{k}}{1 - d^{k-1}}} -
\displaystyle\sum_{l=3}^{k}\displaystyle\frac{\alpha_l}{d}
\ u^{
\dsp\frac{d^{l-k-1}-d^{l-2}-(d-1)}{1-d^{k-1}}}\\
&\ &\\
&+ & s_1 + \cdots + s_{k-1}\\
&\ &\\
&+ & \displaystyle\sum_{j=1}^{\infty} \Big(-\frac{\alpha_2}{d}\Big)^{ j}
u^{
\displaystyle\frac{1-d^k}{1-d^{k-1}} \cdot (d^{1-k})^j} -
\displaystyle\sum_{j=1}^{\infty} \Big(\frac{\Lambda_{1}^{ j} +
\cdots + \Lambda_{k-1}^{ j}}{k-1}\Big) u^{
\displaystyle\frac{1-d^k}{1-d^{k-1}} \cdot d^{-j}}\\
&\ &\\
&- & \displaystyle\sum_{l=3}^{k}\displaystyle\frac{\alpha_l}{d}
\displaystyle\sum_{j=1}^{\infty} \Big(-\frac{\alpha_2}{d}\Big)^{ j} \
u^{
\dsp\frac{d^{l-k-1}-d^{l-2}-(d-1)}{1-d^{k-1}} \cdot (d^{1-k})^j}\\
&\ &\\
&+ & \displaystyle\sum_{l=3}^{k}\displaystyle\frac{\alpha_l}{d}
\displaystyle\sum_{j=1}^{\infty} \Big(\dsp\frac{\Lambda_{1}^{ j} +
\cdots + \Lambda_{k-1}^{ j}}{k-1}\Big) \
u^{
\dsp\frac{d^{l-k-1}-d^{l-2}-(d-1)}{1-d^{k-1}} \cdot d^{-j}}\Bigg]\\
&\ &\\
&=& u^{\displaystyle\frac{(d-1)d^{k-1}}{1 - d^{k-1}}} \
\Bigg[ u^{\displaystyle\frac{1-d^{k}}{1 - d^{k-1}}} -
\displaystyle\sum_{l=3}^{k}\displaystyle\frac{\alpha_l}{d}
u^{
\dsp\frac{d^{l-k-1}-d^{l-2}-(d-1)}{1-d^{k-1}}}\\
&\ &\\
&\ &\qquad \qquad\qquad \quad + \ s_1 \ + \ \cdots\ + \ s_{k-1}\Bigg]\\
&=& u^{\displaystyle\frac{(d-1)d^{k-1}}{1 - d^{k-1}}} \
\Bigg[ s_1 + \dsp\frac{u^{\displaystyle\frac{1-d^{k}}{1 -
d^{k-1}}}}{k-1}-
\displaystyle\sum_{l=3}^{k}\displaystyle\frac{\alpha_l}{d(k-1)}
u^{
\dsp\frac{d^{l-k-1}-d^{l-2}-(d-1)}{1-d^{k-1}}}\\
&\ &\\
&\ & \qquad + \cdots + \ s_{k-1} +
\dsp\frac{u^{\displaystyle\frac{1-d^{k}}{1 - d^{k-1}}}}{k-1} -
\displaystyle\sum_{l=3}^{k}\displaystyle\frac{\alpha_l}{d(k-1)}
u^{
\dsp\frac{d^{l-k-1}-d^{l-2}-(d-1)}{1-d^{k-1}}}
\Bigg]\\
\end{eqnarray*}

Thanks to the change of variable $v := u^{ \dsp\frac{1}{(d^{k-1}-1)\cdot
d^{k-2}}}$,
$\omega$ becomes\\

$$ \omega\ : \big(\mathbb{C} \backslash \overline{\Delta}\big)
\times
\mathbb{C}^{k-1} \ \ \longrightarrow \ \ \big(\mathbb{C} \backslash
\overline{\Delta}\big)
\times
\mathbb{C}^{k-1}$$

\begin{eqnarray*}
\omega(v,s_1,\cdots,s_{k-1}) &=& \Bigg(v^{ d}, \Lambda_1 \Big(s_1+
\dsp\frac{v^{(d^{k}-1)\cdot d^{k-2}}}{k-1}-
\displaystyle\sum_{l=3}^{k}\displaystyle\frac{\alpha_l}{d(k-1)}
v^{\big(-d^{l-3}+d^{l+k-4}+(d-1)d^{k-2}\big)}\Big)\\
&\ &\\
& ,\cdots, & \Lambda_{k-1} \Big(s_{k-1} +
\dsp\frac{v^{(d^{k}-1)\cdot d^{k-2}}}{k-1}-
\displaystyle\sum_{l=3}^{k}\displaystyle\frac{\alpha_l}{d(k-1)}
v^{\big(-d^{l-3}+d^{l+k-4}+(d-1)d^{k-2}\big)} \Big)\Bigg)
\end{eqnarray*}

\smallskip

In variables $(v,s_1,\cdots,s_{k-1})$, the function $f$ renamed $g$
becomes

\begin{eqnarray*}
g(v,s_1,\cdots,s_{k-1}) &=& v^{(d^{k-1}-1)\cdot d^{k-2}} +
v^{(1-d)\cdot d^{2k-3}}\Big(s_1 + \cdots + s_{k-1}\\
&\ & - \displaystyle\sum_{l=3}^{k}\displaystyle\frac{\alpha_l}{d}
v^{-\big(d^{l-3}-d^{l+k-4}-(d-1)d^{k-2}\big)} \Big) + \cdots\\
\end{eqnarray*}

This provides us with the first term

$$g_0 = v^{(d^{k-1}-1)\cdot d^{k-2}} +
v^{(1-d)\cdot d^{2k-3}}\Big(s_1 + \cdots + s_{k-1} -
\displaystyle\sum_{l=3}^{k}\displaystyle\frac{\alpha_l}{d}
v^{\big(-d^{l-3}+d^{l+k-4}+(d-1)d^{k-2}\big)} \Big)$$

of a sequence of functions $(g_n)_n$ which is given by the induction
formula:
$$g_{n+1} = \dsp\frac{1}{\alpha_2} \Big[ g_n\circ\omega^{\circ k} -
\big( g_n\circ\omega^{\circ (k-1)}\big)^d - \alpha_3 g_n\circ \omega -
\cdots - \alpha_k g_n\circ \omega^{\circ (k-2)}\Big].$$

\medskip

In the case $d=2$, by gathering the fractional exponents, we obtain after
calculations

\begin{eqnarray*}
z_{1,n} &=& u^{\displaystyle\frac{2^{k-1}}{1 - 2^{k-1}}} \
\Bigg[ s_1 + \dsp\frac{u^{\displaystyle\frac{1-2^{k}}{1 -
2^{k-1}}}}{k-1}-
\displaystyle\sum_{l=3}^{k}\displaystyle\frac{\alpha_l}{2(k-1)}
u^{
\dsp\frac{2^{l-k-1}-2^{l-2}-1}{1-2^{k-1}}}\\
&+ & \dsp\frac{\alpha_3\alpha_k}{4(k-1)} \ u^{
\displaystyle\frac{-1}{(1-2^{k-1})\cdot 2^{k-1}}} + \cdots + \ s_{k-1} +
\dsp\frac{u^{\displaystyle\frac{1-2^{k}}{1 - 2^{k-1}}}}{k-1}\\
&\ &\\
& -&
\displaystyle\sum_{l=3}^{k}\displaystyle\frac{\alpha_l}{2(k-1)}
u^{
\dsp\frac{2^{l-k-1}-2^{l-2}-1}{1-2^{k-1}}} +
\dsp\frac{\alpha_3\alpha_k}{4(k-1)} \ u^{
\displaystyle\frac{-1}{(1-2^{k-1})\cdot 2^{k-1}}}
\Bigg]\\
\end{eqnarray*}

With the change of variable $v := u^{ \dsp\frac{1}{(2^{k-1}-1)\cdot
d^{k-1}}}$,
$\omega$ becomes\\

$$ \omega\ : \big(\mathbb{C} \backslash \overline{\Delta}\big)
\times
\mathbb{C}^{k-1} \ \ \longrightarrow \ \ \big(\mathbb{C} \backslash
\overline{\Delta}\big)
\times
\mathbb{C}^{k-1}$$

\begin{eqnarray*}
\omega(v,s_1,\cdots,s_{k-1}) &=& \Bigg(v^{ 2}, \Lambda_1 \Big(s_1+
\dsp\frac{v^{(2^{k}-1)\cdot 2^{k-1}}}{k-1}-
\displaystyle\sum_{l=3}^{k}\displaystyle\frac{\alpha_l}{2(k-1)}
v^{\big(-2^{l-2}+2^{l+k-3}+2^{k-1}\big)}\\
&\ &\\
&+& \displaystyle\frac{\alpha_3\alpha_k}{4(k-1)} v\Big)\ ,\cdots,
\Lambda_{k-1} \Big(s_{k-1} +
\dsp\frac{v^{(2^{k}-1)\cdot 2^{k-1}}}{k-1}\\
&\ &\\
&-&
\displaystyle\sum_{l=3}^{k}\displaystyle\frac{\alpha_l}{2(k-1)}
v^{\big(-2^{l-2}+2^{l+k-3}+2^{k-1}\big)} +
\displaystyle\frac{\alpha_3\alpha_k}{4(k-1)} v \Big)\Bigg)
\end{eqnarray*}

\smallskip

In variables $(v,s_1,\cdots,s_{k-1})$, the function $f$ renammed $g$
becomes in the case $d=2$

\begin{eqnarray*}
g(v,s_1,\cdots,s_{k-1}) &=& v^{(2^{k-1}-1)\cdot 2^{k-1}} +
v^{ -2^{2k-2}}\Big(s_1 + \cdots + s_{k-1}\\
&- & \displaystyle\sum_{l=3}^{k}\displaystyle\frac{\alpha_l}{2}
v^{-\big(2^{l-2}-2^{l+k-3}-2^{k-1}\big)} +
\frac{\alpha_3\alpha_{k}}{4} v\Big) +
\cdots\\
\end{eqnarray*}

This provides us with the first term

$$g_0 = v^{(2^{k-1}-1)\cdot 2^{k-1}} +
v^{- d^{2k-2}}\Big(s_1 + \cdots + s_{k-1} -
\displaystyle\sum_{l=3}^{k}\displaystyle\frac{\alpha_l}{2}
v^{\big(-2^{l-2}+2^{l+k-3}+2^{k-1}\big)} + \frac{\alpha_3\alpha_{k}}{4}
v\Big)$$

of a sequence of functions $(g_n)_n$ given by the induction formula:
$$g_{n+1} = \dsp\frac{1}{\alpha_2} \Big[ g_n\circ\omega^{\circ k} -
\big( g_n\circ\omega^{\circ (k-1)}\big)^2 - \alpha_3 g_n\circ \omega -
\cdots - \alpha_k g_n\circ \omega^{\circ (k-2)}\Big].$$

The sequence $(g_n)_n$ converges uniformally
towards the function $g$ on a domain of the kind:\\

$$ \left\{ \begin{array}{l}
\vert v\vert \geq 1+\varepsilon\\
\ \\
\vert s_1\vert \leq K_1\vert v\vert^{(d^{k-1}-1)\cdot d^{k-2}} \\
\ \vdots\\
\ \vdots\\
\vert s_{k-1}\vert \leq K_{k-1}\vert v\vert^{(d^{k-1}-1)\cdot d^{k-2}} \\

\end{array}\right.$$
or

$$ \left\{ \begin{array}{l}
\vert v\vert \geq 1+\varepsilon\\
\ \\
\vert s_1\vert \leq K_1\vert v\vert^{(2^{k-1}-1)\cdot 2^{k-1}} \\
\ \vdots\\
\ \vdots\\
\vert s_{k-1}\vert \leq K_{k-1}\vert v\vert^{(2^{k-1}-1)\cdot 2^{k-1}} \\

\end{array}\right.$$

\medskip

Let us set

$$G \ = \left[ \begin{array}{l}
\quad g\ \circ\ \omega^{\circ (k-1)}\\
\\
\quad g \\
\\
\quad g\ \circ\ \omega\\
\\
\vdots
\\
\quad g\ \circ\ \omega^{\circ (k-2)}
\end{array}
\right]$$\\

then $G$ is such that: $G\ \circ\ \omega = H\ \circ \ G$.\\

Let us give now the proof of the proposition \ref{P}. We only give it in
the case
$d\geq 3$. The proof would be almost the same for $d=2$.\\

\demo Let $(z_{1,n})_n$ and $(\zeta_{1,n})_n$ be two sequences which
tend to infinity and are such that:
$$\left\{ \begin{array}{l}
z_{1,n+1} = z_{1,n}^d + \alpha_2\ z_{1,n-k+1} + \ \cdots + \alpha_k \
z_{1,n-1}\\
\ \\
\zeta_{1,n+1} = \zeta_{1,n}^{d} + \alpha_2\ \zeta_{1,n-k+1} + \ \cdots
+
\alpha_k
\ \zeta_{1,n-1}\\

\end{array}\right.$$

then either $z_{1,n}$ and $\zeta_{1,n}$ do not have the same
$\varphi$ and in this case they have nothing to do one another, or
they have the same $\varphi$ and in this case, $$z_{1,n} - \zeta_{1,n}
\sim {\rm cste}\ \Big(z_{1,n}^{\dsp\frac{(d-1)d^{k-1}}{1-d^{k-1}}}\Big)
\sim {\rm cste}\
\Big(\zeta_{1,n}^{\dsp\frac{(d-1)d^{k-1}}{1-d^{k-1}}}\Big).$$

We are led to set the following definitions:\\

\begin{definition} Let $(z_{1,n})_n$ be a sequence such that $z_{1,n}$
tends to infinity with $n$ and the quotient
$\displaystyle\frac{z_{1,n+1}}{z_{1,n}} \to + \infty$.\\

$i$. $(z_{1,n})$ is said to be class $(k,d)$ if for any $ n$ in
$\mathbb{Z}$, it satisfies
$$ - z_{1,n+1} + z_{1,n}^d + \alpha_2\ z_{1,n-k+1} + \ \cdots + \alpha_k \
z_{1,n-1} = 0.$$

$ii.$ $(z_{1,n})$ is said to be almost class $(k,d)$ if for any $n$,
it satisfies:

$$ - z_{1,n+1} + z_{1,n}^d + \alpha_2\ z_{1,n-k+1} + \ \cdots + \alpha_k \
z_{1,n-1}= {\rm o
o}\Big(z_{1,n}^{\dsp\frac{1}{1-d^{k-1}}}\Big)$$ where
$u_n :=
{\rm o o}(v_n) \Longleftrightarrow \displaystyle\sum_{n_0}^{\infty}\
\left\vert\frac{u_n}{v_n}\right\vert <\ +\ \infty$ for $n_0 >> 1$.
\end{definition}
\smallskip

\begin{definition} Let $(z_{1,n})_n$ and $(\zeta_{1,n})_n$ be almost class
$(k,d)$. The two sequences are:

$\qquad \bullet$ neighbour if $z_{1,n} - \zeta_{1,n} = O\
\Big(z_{1,n}^{\dsp\frac{(d-1)d^{k-1}}{1-d^{k-1}}}\Big),$\\

$\qquad \bullet$ twin if $z_{1,n} - \zeta_{1,n} = o\
\Big(z_{1,n}^{\dsp\frac{(d-1)d^{k-1}}{1-d^{k-1}}}\Big).$\\

These are equivalence relations.
\end{definition}

\begin{Prop} Let $(z_{1,n})$ be an almost class $(k,d)$ sequence. Then
there exists one and only one sequence $(y_{1,n})_n$ class $(k,d)$
twin with respect to $(z_{1,n})_n$
\end{Prop}

We introduce the distance $j_{n_0} (z_{1,n},\ \zeta_{1,n}) =
{\rm sup}_{\ n\geq \ n_0} \vert z_{1,n} -
\zeta_{1,n}\vert \ z_{1,n}^{\dsp\frac{(d-1)d^{k-1}}{d^{k-1}-1}}$.\\

The symmetry follows from the fact we have an equivalence relation.\\ The
distance equals $0$ only when $(z_{1,n})_n$ and
$(\zeta_{1,n})_n$ are equal from a certain rank but it does not satisfy
the triangular inequality because of the non-linearity
with respect to $z_{1,n}^{\dsp\frac{(d-1)d^{k-1}}{d^{k-1}-1}}$.\\

Let $Z = (z_{1,n})$ be an almost class $(k,d)$ sequence. We
associate $\mathcal{F}(Z) = \zeta_{1,n}$ where \\
\begin{center}$\zeta_{1,n} = \displaystyle\frac{1}{\alpha_2} \Big(
z_{1,n+k} - z_{1,n+k-1}^{ d} - \alpha_3 \ z_{1,n+1} - \cdots - \alpha_k
z_{1,n-k+2}
\Big)$ ;
\
$Z_p =
\mathcal{F}^{
\circ p}(Z)$.\end{center}
The fixed points of $\mathcal{F}$ are class $(k,d)$ sequences.
$Z_p$ are almost class $(k,d)$ and twin. There remains to see that
$(Z_p)$ converges uniformally on a domain of the kind $[n_0\ ,\
\infty]$.
$$j_{n_0} (Z_p ,\ Z_{p+\nu}) \ \lesssim
\displaystyle\sum_{l=p}^{p+\nu-1}\ z_{1,l}^{\frac{1}{d^{k-1}-1}}\
\big(\ z_{1,l}^d + \alpha_2 z_{1,l-k+1} + \cdots + \alpha_k\ z_{1,l-1} -
z_{1,l+1} \big) <
\infty$$

hence $(Z_{p})$ is a Cauchy sequence for $j_{n_0}$ so it converges
uniformally for $n\geq n_0$ towards a class $(k,d)$ sequence for
$n\geq n_0$. \\

{\bf Surjectivity of $G$.}\\

If $(z_{1,n})$ is a class $(k,d)$ sequence then there exists $(
v ,\ s_1 ,\ \cdots, \ s_{k-1} ) \in \mathbb{C} \backslash
\overline{\Delta}
\times
\mathbb{C}^{ k-1}$ such that $\mu_n = g_0 \circ \omega^{\circ n} (v,\ s_1
,\ \cdots, \ s_{k-1} ) $ is twin with respect to $(z_{1,n})$. \\

We can suppose

$$\left(\begin{array}{l} z_{1, 0}\\
z_{1,-1}\\
\vdots\\
z_{1,-k+1}\\
\end{array}\right) \ \in V^{ +} \ {\rm and\ take\ }
\varphi\left(\begin{array}{l} z_{1, 0}\\
z_{1,-1}\\
\vdots\\
z_{1,-k+1}\\
\end{array}\right) = v^{(d^{k-1}-1)\cdot d^{k-2}}$$

\begin{eqnarray*}
\omega(v,0,\cdots,0) &=& \Bigg(v^{ d}, \Lambda_1 \Big(
\dsp\frac{v^{(d^{k}-1)\cdot d^{k-2}}}{k-1}-
\displaystyle\sum_{l=3}^{k}\displaystyle\frac{\alpha_l}{d(k-1)}
v^{\big(-d^{l-3}+d^{l+k-4}+(d-1)d^{k-2}\big)}\Big)\\
&\ &\\
& ,\cdots, & \Lambda_{k-1} \Big(
\dsp\frac{v^{(d^{k}-1)\cdot d^{k-2}}}{k-1}-
\displaystyle\sum_{l=3}^{k}\displaystyle\frac{\alpha_l}{d(k-1)}
v^{\big(-d^{l-3}+d^{l+k-4}+(d-1)d^{k-2}\big)} \Big)\Bigg)
\end{eqnarray*}

and by iterating, we obtain

\beqan
\omega^{\circ n}(v,\ 0,\ \cdots,\ 0) &=& \Bigg(v^{ d^{ n}},
\ \frac{1}{k-1}\displaystyle\sum_{j=1}^{n} \Lambda_1^{ j}
\Big(v^{(d^{k}-1)\cdot d^{k-2}\cdot d^{n-j}}\\
&\ & \qquad \qquad \qquad -
\displaystyle\sum_{l=3}^{k}\displaystyle\frac{\alpha_l}{d}
v^{\big(-d^{l-3}+d^{l+k-4}+(d-1)d^{k-2}\big)\cdot d^{n-j}}\Big)\\
&\ & \qquad ,
\cdots, \frac{1}{k-1}\displaystyle\sum_{j=1}^{n} \Lambda_{k-1}^{ j}
\Big(v^{(d^{k}-1)\cdot d^{k-2}\cdot d^{n-j}}\\
&\ & \qquad\qquad\qquad -
\displaystyle\sum_{l=3}^{k}\displaystyle\frac{\alpha_l}{d}
v^{\big(-d^{l-3}+d^{l+k-4}+(d-1)d^{k-2}\big)\cdot d^{n-j}}\Big)\Bigg)
\eeqan

Let us set
$$\chi_n := g_{ 0} \circ \omega^{\circ n}(v,\ 0,\ \cdots,\ 0)$$

then, $(z_{1,n} - \chi_n)\cdot\chi_n^{(d-1)\cdot d^{2k-3}}$ converges
to $w$.\\

Let now $\mu_n$ be $\mu_n := g_{ 0} \circ \omega^{\circ n}(v,\ w,\
\cdots,\ w)$.\\

The sequences $(z_{1,n})$ and $(\mu_n)$ are twin because $$z_{1,n} -
\mu_n = O\
\Big(\mu_n^{\dsp\frac{(d-1)d^{k-1}}{1-d^{k-1}}}\Big)$$

hence $G(u,\ w,\ \cdots,\ w) = \big( z_{1, 0},
z_{1,-1},
\cdots,
z_{1,-k+1}\big)$. The map $G$ is surjective.\\

\medskip

As for the fibers: if $(v,\ s_1,\ \cdots,\ s_{k-1})$ and $(v',\
s'_1,\ \cdots,\ s'_{k-1})$ are such that for any $ n\geq 0\ ,\
\omega^{\circ n}(v,\ s_1,\ \cdots,\ s_{k-1}))
\neq
\omega^{\circ n} (v',\ s'_1,\ \cdots,\ s'_{k-1})$then the sequences $g_{
0}
\circ
\omega^{\circ n}(v,\ s_1,\ \cdots,\ s_{k-1})$ and $g_{
0} \circ \omega^{\circ n}(v',\ s'_1,\ \cdots,\ s'_{k-1})$ are not twin
and consequently

$$G(v, \ s_1,\ \cdots,\ s_{k-1})
\neq G(v',\ s'_1,\ \cdots,\ s'_{k-1}). $$

hence the fiber of $G(v, \ s_1,\ \cdots,\ s_{k-1})$ is made of
k-uplets $(v', \ s'_1,\ \cdots,\ s'_{k-1})$ where
$\displaystyle\frac{v'}{v}$ is a $d^{n}$-th root of unity
and where for $i$ from $1$ to $k-1$:

\beqan
&\ & \Lambda_i^{ n} s_i +
\displaystyle\sum_{j=1}^{n}\dsp\frac{\Lambda_i^{ j}}{k-1}\
\Big(v^{(d^{k}-1)\cdot d^{k-2}\cdot d^{n-j}} -
\displaystyle\sum_{l=3}^{k}\displaystyle\frac{\alpha_l}{d}
v^{\big(-d^{l-3}+d^{l+k-4}+(d-1)d^{k-2}\big)\cdot d^{n-j}}\Big)\\
&=& \quad \Lambda_i^{ n} s'_i +
\displaystyle\sum_{j=1}^{n}\dsp\frac{\Lambda_i^{ j}}{k-1}\
\Big(v^{'(d^{k}-1)\cdot d^{k-2}\cdot d^{n-j}} -
\displaystyle\sum_{l=3}^{k}\displaystyle\frac{\alpha_l}{d}
v^{'\big(-d^{l-3}+d^{l+k-4}+(d-1)d^{k-2}\big)\cdot d^{n-j}}\Big)\\
\eeqan

or as well

\beqan
\Delta_i\ (v,\ v') &:=& s_i' - s_i\\
&=&\displaystyle\sum_{j=1}^{n}\dsp\frac{\Lambda_i^{ j-n}}{k-1}\
\Bigg(v^{(d^{k}-1)\cdot d^{k-2}\cdot d^{n-j}} - v^{'(d^{k}-1)\cdot
d^{k-2}\cdot d^{n-j}}\\
&\ & \qquad -\displaystyle\sum_{l=3}^{k}\displaystyle\frac{\alpha_l}{d}
v^{\big(-d^{l-3}+d^{l+k-4}+(d-1)d^{k-2}\big)\cdot d^{n-j}}\\
&\ & \qquad +
\displaystyle\sum_{l=3}^{k}\displaystyle\frac{\alpha_l}{d}
v^{'\big(-d^{l-3}+d^{l+k-4}+(d-1)d^{k-2}\big)\cdot d^{n-j}}\Bigg)\\
& =&
\sum_{m=0}^{n-1}\ \dsp\frac{\Lambda_i^{ -m}}{k-1}\
\Bigg(v^{(d^{k}-1)\cdot d^{k-2}\cdot d^{m}} - v^{'(d^{k}-1)\cdot
d^{k-2}\cdot d^{m}}\\
&\ & \qquad - \displaystyle\sum_{l=3}^{k}\displaystyle\frac{\alpha_l}{d}
v^{\big(-d^{l-3}+d^{l+k-4}+(d-1)d^{k-2}\big)\cdot d^{m}}\\
&\ & \qquad +
\displaystyle\sum_{l=3}^{k}\displaystyle\frac{\alpha_l}{d}
v^{'\big(-d^{l-3}+d^{l+k-4}+(d-1)d^{k-2}\big)\cdot d^{m}}\Bigg)\\
\eeqan

finally, the fiber is countable.\\

What does it look like? The term $v'$ describes a dense part of
the cercle with a radius equal to $\vert v\vert$. The fiber is
discrete. Indeed, when
$v' = e^{\ 2i\pi p\cdot d^{ -n}}\cdot v$ with a very big $n$,
$\Delta_i\ (v,\ v')$ is very big (if $p$ is such that $2p\cdot d^{
-n}\cdot(d^{k}-1)\cdot d^{k-2}\cdot d^{n-1}$ is odd).
\begin{eqnarray*}
\Delta_i\ (v,\ v')
&\stackrel{n \ big,p\
chosen}{\sim}& \dsp\frac{\Lambda_i^{ 1-n}}{k-1}\
\Bigg(v^{(d^{k}-1)\cdot d^{k-2}\cdot d^{n-1}} - v^{'(d^{k}-1)\cdot
d^{k-2}\cdot d^{n-1}}\Bigg)\\
&\ &\\
&\sim& \dsp\frac{2\Lambda_i^{ 1-n}}{k-1}\
\cdot v^{(d^{k}-1)\cdot d^{k-2}\cdot d^{n-1}}
\end{eqnarray*}

\medskip

The action of $\displaystyle\mathbb{Z}\ \Big[\ \frac{1}{d}\ \Big]\ /
\mathbb{Z}$ defined by:

$$\begin{array}{cccc} \mathcal{F} : &\left[\displaystyle\mathbb{Z}\
\Big[\
\frac{1}{d}\ \Big]\ / \mathbb{Z}\ \right]\ \times\ \Big[\big(\mathbb{C}\
\backslash\overline{\Delta}\big)\ \times\ \mathbb{C}^{ k-1}\Big]
&\longrightarrow & \big(\mathbb{C}\
\backslash\overline{\Delta}\big)\ \times\ \mathbb{C}^{ k-1}\\
&\ &\ &\ \\
&(\ \theta , \ v ,\ s_1 ,\ \cdots,\ s_{k-1} )
&\longmapsto &\mathcal{F}\ (\ \theta, \ v , \ s_1 ,\ \cdots,\ s_{k-1})\\
&\ &\ & =
(v' ,\ s'_1 ,\ \cdots,\ s'_{k-1})
\end{array}$$
with $\left\{\begin{array}{l}
v' = v\cdot e^{ 2i\pi \theta}\\
\\
s'_i = s_i\ + \ \Delta_i\ (v,\ v')\\
\end{array}\right.$
whose orbits are the fibers of $G$, is free and effective. The fibers
are homeomorphic to $\displaystyle\mathbb{Z}\ \Big[\
\frac{1}{d}\ \Big]\ / \mathbb{Z}$.\\

The bundle is locally trivial as a consequence of $\varphi$ being
locally defined. We will deal later in the appendix the problems of
convergence.\qed

\medskip
Finally, we have established the following theorem:

\begin{theorem}\label{pro}
$$\pi_1\ ( U^{+}) \simeq \mathbb{Z}\ \Big[\ \frac{1}{d}\ \Big]$$

$$\widetilde{U^{ +}} = \mathbb{H}\ \times \mathbb{C}^{ k-1}$$

where $\mathbb{H}$ denotes the Poincar\'e half-plane.
\end{theorem}

\begin{remark}

The attraction basin $U^+$ is a domain of holomorphy. Moreover,
thanks to the automorphisms obtained from the action of $\pi_1(U^{+})$ on
$\mathbb{H} \times \mathbb{C}^{ k-1}$, we note that the action
of $\ \mathbb{C}^{k-1}$ on $\mathbb{H} \times \mathbb{C}^{ k-1} $ induces
an action of $\ \mathbb{C}^{k-1}$ on $U^+$. The orbits of this action
give a foliation of codimension $1$ on $U^+$.
\end{remark}

\medskip

\subsection{Appendix}$\ $\\

In this part, we want to solve the problems of convergence which have
appeared in our proofs. In the proof of the lemma \ref{L}, we said that
the sequence of functions defined by $g_0$ and the induction relation
given by:

$$g_{n+1} = \dsp\frac{1}{\alpha_2} \Big[ g_n\circ\omega^{\circ k} -
\big( g_n\circ\omega^{\circ (k-1)}\big)^d - \alpha_3 g_n\circ \omega -
\cdots - \alpha_k g_n\circ \omega^{\circ (k-2)}\Big].$$

was formally convergent. We have to explicit the notion of convergence.\\

Let $K$ be a field or an integral ring. We set $K \lceil X \rceil$
the space of the formal linear combinations $\displaystyle\sum_{\kappa\in
A} \ k_{\kappa} \ X^{\kappa}$ where $A$ is a closed part of $\mathbb{R}$,
lowerly bounded and discrete on the right that is to say
$$\forall
\
\kappa \ \in \ A\ ,\ \exists\ \varepsilon > 0
\ , \
]\ \kappa\ ,\ \kappa + \varepsilon\ [\ \cap \ A = \emptyset.$$
We can define the sum and the product of two elements of $K \lceil X
\rceil$ providing this latter with a structure of
$K$-algebra. We could also define $K \lfloor X \rfloor$
by identifying $K \lceil X^{ -1} \rceil$ and $K \lfloor X
\rfloor$.

\begin{center} Let $f$ be  in $K \lceil X \rceil$, $f \neq 0$.
\end{center}

Let us note $f = \displaystyle\sum_{\kappa\in\mathbb{R}} \ k_{\kappa} \
X^{\kappa}$ then the number
min$\{ \kappa\in\mathbb{R} \ ,\ k_{\kappa}\neq 0\}$ exists and is finite.
We call it $v_X\ (f)$ and we set $v_X \ (0) = +\infty$.\\ By extending
this notion,
$v_{X^{\kappa}}(f) = \frac{1}{\kappa} \ v_X (f)$. $v_X$ has the
properties of a valuation. We can define a norm on $K \lceil X \rceil$
given by: $\vert\vert f \vert\vert
= exp \ [ \ -v_X \ (f) \ ]$; it is a norm because

$$\vert\vert f \vert\vert
= 0 \Leftrightarrow -v_X \ (f) \ = -\infty \Leftrightarrow v_X \ (f) \ =
\infty
\Leftrightarrow f = 0$$ as for the triangular inequality, it is a
consequence of the convexity of the function exponential.
$K \lceil X \rceil$ with this norm is complete.\\ Moreover, we have
the equality: $v_X \ (f g) = v_X \ (f) + v_X \ (g)$ as a consequence of
$K$ being integral. Finally, $K \lceil X \rceil$ is integral.

\begin{eqnarray*}
\vert\vert fg \vert\vert
&=& exp \ [ \ -v_X \ (fg) \ ] = exp \ [ \ -v_X \ (f) \ - v_X \ (g) ]\\
&=&
exp \ [ \ -v_X \ (f)\ ] \ exp \ [ \ -v_X \ (g)\ ] = \vert\vert f
\vert\vert \ \vert\vert g \vert\vert
\end{eqnarray*}

If $K$ is a field, so is $K \lceil X \rceil$ and $x\ \mapsto
x^{ -1}$ is continuous.\\

In the space $\mathbb{C}[s_1,\cdots,s_{k-1}]\lceil u^{ -1} \ \rceil =
\mathbb{C}[s_1,\cdots,s_{k-1}]\lfloor u \ \rfloor$, the formula
which defines $g_0$ takes sense and the functions $g_n$ are defined by
the induction formula.\\ We show that
$$g_1 = g_0 + o \ \Big(u^{ \frac{(d-1)d^{k-1}}{1-d^{k-1}} }\Big)$$
and consequently
$$g_1 = g_0 + O \
\Big(u^{ \frac{(d-1)d^{k-1}}{1-d^{k-1}} + \ \varepsilon \ }\Big).$$ We
transform the sequence into a series by setting $\Delta_k = g_k -
g_{k-1}$. Then, we have:

\begin{eqnarray*} \Delta_{k+1} &=& g_{k+1} - g_{k}\\
&=& \ \ \frac{1}{\alpha_2}\Big[
g_k(u^{d^k}) - \big(g_k(u^{d^{k-1}})\big)^d - \alpha_3 g_{k}(u^d) -
\cdots - \alpha_k g_k\big(u^{d^{k-2}}\big)\
\Big]\\
&\ & - \frac{1}{\alpha_2}\Big[
g_{k-1}(u^{d^k}) - \big(g_{k-1}(u^{d^{k-1}})\big)^d - \alpha_3
g_{k-1}(u^d) -
\cdots - \alpha_k g_{k-1}\big(u^{d^{k-2}}\big)\
\Big]\\
&=& \frac{1}{\alpha_2}\Big[
\Delta_k(u^{d^k}) - \Delta_k(u^{d^{k-1}})\big( g_k^{d-1} + \cdots +
g_{k-1}^{d-1}\big) -
\alpha_3
\Delta_{k}(u^d) -
\cdots - \alpha_k \Delta_k\big(u^{d^{k-2}}\big)\
\Big].\\
\end{eqnarray*}

We prove by induction that $v_{u^{ -1}} (\Delta_k) \geq
\frac{(d-1)d^{k-1}}{d^{k-1}-1} + d^{ k} \varepsilon$.\\

One has $v_{u^{ -1}} (g_k^{d-1} + \cdots + g_{k-1}^{d-1}) = 1-d$.

$$v_{u^{-1}} (\Delta_k(u^{d^k})) =
d^k\Big(\frac{(d-1)d^{k-1}}{d^{k-1}-1} + d^{ k} \varepsilon\Big) >
\frac{(d-1)d^{k-1}}{d^{k-1}-1} + d^{ k+1} \varepsilon.$$

\begin{eqnarray*}v_{u^{ -1}} \Big(\big( g_k^{d-1} + \cdots +
g_{k-1}^{d-1}\big)\Delta_{k}(u^{d^{k-1}})\Big) &=& v_{u^{
-1}}\big(g_k^{d-1} + \cdots + g_{k-1}^{d-1}\big) + v_{u^{ -1}}
\Big(\Delta_k(u^{d^{k-1}})\Big)\\ &\ &\\
&>& \frac{(d-1)d^{k-1}}{d^{k-1}-1} + d^{ k+1} \varepsilon.
\end{eqnarray*}

$$v_{u^{-1}} (\Delta_k(u^d)) =
d\Big( \frac{(d-1)d^{k-1}}{d^{k-1}-1} + d^{ k} \varepsilon\Big)
 	 	\frac{(d-1)d^{k-1}}{d^{k-1}-1} + d^{ k+1} \varepsilon.$$ The
series $\Delta_k$ is convergent because the general term tends to $0$
thus so does the sequence of $(g_k)$. It converges towards
an element of $\mathbb{C}[s_1,\cdots,s_{k-1}] \lceil u^{-1} \rceil$.\qed\\

\section{The Case of dimension $3$.}

\smallskip

Now, we deal with the particular case when the dimension $k$ equals $3$
and the degree $d$ equals $2$. The quadratic
automorphisms have been classified by Fornaess and Wu
\cite{FW}. We are
able to list all the automorphisms of this classification for which the
method described before applies.

\begin{theorem}
The quadratic automorphisms of $\mathbb{C}^3$ in the list of
Fornaess-Wu \cite{FW} for which the method described previously
works are:

\begin{enumerate}
\item in the first class, those of the form $$ H_1(x,\ y,\ z) = \left\{
\begin{array}{l}
\alpha x^2 + bx +ay + \gamma \\
\\
x + \varepsilon \\
\\
\nu z\\
\end{array}\right.$$

with $ a\neq 0,\ \mid a\mid <2,\ \nu
\neq 0,\ \alpha \neq 0$ and $b,\ \gamma,\ \varepsilon$
any constant.\\
For this class, there are also those of the form

$$\qquad H_1(x,\ y,\ z) = \left\{
\begin{array}{l}
\alpha x^2 +ay + p_1(z) x + p_2(z)\\
\\
\varepsilon + x \ \ \qquad \qquad \qquad \qquad,\ a\neq 0,\ \mid a\mid
<2,\ \alpha \neq 0\ \ \ \\
\\
z\\
\end{array}\right.$$
$p_1$ and $p_2$ are polynomes with degrees deg$(p_1) \leq 1$,
deg$(p_1) \leq 2$.\\

\item In the second class, those of the form
$$ \quad H_2(x,\ y,\ z) = \left\{
\begin{array}{l}
a x + P(y,\ z)\\
\\
\alpha y^2 + \beta y + b z + c\ \ \qquad , \ a\neq 0,\ \alpha\neq 0 ,\
b\neq 0,\
\mid b\mid <2\\
\\
y\\
\end{array}\right.$$
$P$ is a polynome with a degree smaller than $2$,  $\beta$
and $c$ are any constant.\\

\item In the third class, those of the form
$$ \quad \qquad H_3(x,\ y,\ z) = \left\{
\begin{array}{l}
\alpha x^2 + \mu x + \nu z + ay + \delta\\
\\
\varepsilon x + z + \rho \ \ \qquad \qquad \quad ,\ a\neq 0,\ \mid
a\mid <2,\ \alpha \neq 0\\
\\
x\\
\end{array}\right.$$\\
.

\item In the fourth class, those of the form

$$ \ H_4(x,\ y,\ z) = \left\{
\begin{array}{l}
\alpha x^2 + \gamma x + \gamma' y + az + \delta \\
\\
x + \rho \ \ \qquad \qquad \qquad ,\ a\neq 0,\ \mid a\mid <2,\
\alpha \neq 0\\
\\
y\\
\end{array}\right.$$
there is no condition on the other constants.\\
There are also those of the form:

$$ \qquad \quad \ H_4(x,\ y,\ z) = \left\{
\begin{array}{l}
\gamma y + az + \varepsilon \\
\\
\alpha y^2 +\ \nu y + \ x + \delta\ \qquad \quad ,\ a\neq 0,\
\mid a\mid <2,\ \alpha \neq 0\\
\\
y\\
\end{array}\right.$$\\
There is no condition on the other constants.

\item In the fifth class, those of the form

$$ \qquad \qquad H_5(x,\ y,\ z) = \left\{
\begin{array}{l}
\alpha x^2 + \nu x + \delta \ + az\\
\\
\beta x^2 + \gamma x + \rho \ + b y \ \qquad \quad ,\ a\ {\rm et} \ b\neq
0,\ \mid a\mid <2,\ \alpha \neq 0\\
\\
x\\
\end{array}\right.$$
There is no condition on the other constants.
\end{enumerate}

Each of these automorphisms admits an attracting fixed point at infinity
whose attracting basin $U^+$ is such that:\\
$${\it the\  fundamental \ group \ }\ \pi_1(U^+) \simeq
\mathbb{Z} \ \Big[\frac{1}{2}\Big]$$

\begin{center}
{\it and\ his\ covering } $\widetilde{U^+} \simeq \mathbb{H}\times
\mathbb{C}^2$
\end{center}

%\\

\qquad \qquad \qquad \qquad where $\mathbb{H}$ denotes the Poincar\'e
half-plane.
\end{theorem}

\medskip

\demo Let us set, for $n \in \mathbb{Z}$, $(x_n,\ y_n,\ z_n) := H^{\circ
n} (x,\ y,\ z)$. In order to study the sequence of the iterates of $(x,\
y,\ z)$ by $H$, all we need is to study one of the sequences $(x_n)_n,\
(y_n)_n$ or $(z_n)_n$.  If the induction relation that gives the $n+1$-
indexed term only depends on the $n$th et
$n-1$th ones, we will say that {\it the automorphism can be considered
as an H\'enon}. Otherly, the induction formula depends on the terms of
indexes $n$, $n-1$ and $n-2$.

We distinguish two cases. The automorphisms that can be considered
as H\'enon ones and that can be dealt with, in the same way as in
T.Bousch's thesis \cite{B};  the others that can be dealt with in the
same way as what precedes.

As far as the automorphisms that can be considered as H\'enon are
concerned, there are, in the previous list, those of the {\bf $1$st,
$2$nd and $5$th classes}. Those of the {\bf  $3$rd and $4$th classes}
belong to the second category.\qed

\medskip

\begin{remark} We have to notice that those automorphisms satisfy the
property \ref{P}. In particular, $\omega$ is given explicitely.

\end{remark}

\medskip

We end up with a question:\\
\smallskip

{\bf Question :} Let $f$ be $f : \mathbb{C}^k \longrightarrow
\mathbb{C}^k$  a polynomial automorphism of $\mathbb{C}^k$ which admits a
fixed attracting point at infinity and let $U^+$ be its attracting basin,
do we always have:

$$ \widetilde{U^+} \simeq \mathbb{H}\times
\mathbb{C}^{k-1} ?$$

\end{document}